# Two-Stage Stochastic Optimization Frameworks to Aid in Decision-Making Under Uncertainty for Variable Resource Generators Participating in a Sequential Energy Market

Razan A. H. Al-Lawati[1], Jose L. Crespo-Vazquez[2], Tasnim Ibn Faiz[1], Xin Fang[3], Md. Noor-E-Alam[1*]

[1]Dept. of Mechanical and Industrial Engineering, Northeastern University, 360 Huntington Avenue, 02115, Boston, MA, USA

[2]Copernicus Institute of Sustainable Development, Utrecht University, Vening Meinesz building, Princetonlaan 8a, 3584 CB Utrecht, Netherlands

[3]National Renewable Energy Laboratory (NREL), 15013 Denver West Parkway, Golden, CO 80401

* Corresponding author Tel.: +1-617-373-2275

E-mail address: mnalam@neu.edu

Abstract

Decisions for a variable renewable resource generator's commitment in the energy market are typically made in advance when little information is obtainable about wind availability and market prices. Much research has been published recommending various frameworks for addressing this issue. However, these frameworks are limited as they do not consider all markets a producer can participate in. Moreover, current stochastic programming models do not allow for uncertainty data to be updated as more accurate information becomes available. This work proposes two decision-making frameworks for a wind energy generator participating in day-ahead, intraday, reserve, and balancing markets. The first framework is a two-stage stochastic convex optimization approach, where both scenario-independent and scenario-dependent decisions are made concurrently. The second framework is a series of four two-stage stochastic optimization models wherein the results from each model feed into each subsequent model allowing for scenarios to be updated as more information becomes available to the decision-maker. In the simulation experiments, the multi-phase framework performs better than the single-phase in every run, and results in an average profit increase of 7%. The proposed optimization frameworks aid in better decision-making while addressing uncertainty related to variable resource generators and maximize the return on investment.

## 1. Introduction

Renewable energy has started to play a more prominent role in the electricity markets [1]; however, additional complexity arises when considering Variable Renewable Resource Generators (VRRGs) due to the variable nature of the available resource. Wind farms and solar photovoltaic plants are considered VRRGs because unlike other renewable resources, such as biogas and geothermal, the amount of the resource available is uncertain and dependent on nature. It cannot be controlled and is difficult to predict. Moreover, most electricity markets around the world demand that the amount of energy sold by a generator be committed ahead of time of actual participation, typically, one day before the actual buying and selling of that energy [2]. At the time of decision making, information about wind availability and market prices is uncertain. The energy producers are subject to monetary penalties if they deviate from their committed schedules [3]. For this reason, decision-making becomes a complex task.

Most electricity markets are designed sequentially to reduce discrepancies between energy committed and energy available when more information is made available [4],[5]. However, as the time gets closer to the point of actual participation, market prices become less favorable to the generator. Hence, the energy producer must create accurate forecasts of energy availability ahead of time [2], in addition to following the right strategy for participating in suitable markets to make the most profit [6].

To make the work presented in this paper as generalizable and comprehensive as possible, we consider a typical European market that contains all of the markets an energy producer can participate in. Although electricity markets around the world can vary, they all contain the same general market structures [7]. The pool market (PM) consists of three different markets in which commitments can be sequentially updated: (1) a day-ahead market (DAM), (2) adjustment markets, such as an intraday market (IDM), and (3) a balancing market (BM) [8]. Additionally, separately from the PM, there are other markets such as the reserve and regulation markets to ensure secure system operation and energy delivery [8]. In this paper, we consider the regulation market wherein up and down real-time load-following capability is provided to enforce a continuous balance between production and consumption [8]. The regulation market (RM) is typically cleared once a day on an hourly basis and assigns production units the power bands to be used in real-time operation for load following [8].

The electricity market is usually regulated by a System Operator (SO), which has the role of maintaining the reliability, security, stability, and quality of the power supply to the customers [9]. Matching the supply and demand is vital for electrical grid operation and is termed "frequency regulation". One method of frequency generation is requiring the generator to be flexible and increase or decrease output by some amount [10]. This is known as "regulation up" and "regulation down". The wind farm (WF) is subject to the cost of regulation if the generated energy in real-time is less than the committed volume when the up-regulation price is activated or more than the committed volume when the down-regulation price is activated [3].

Another way of overcoming the problem of uncertainty is by using Energy Storage Systems (ESS). Excess energy can be stored and used later when natural resources are not available. However, ESS are expensive and do not always make up for the additional profit that the generator makes. Research has shown the theoretical and practical significance of integrating storage to a WF, mainly for adding an arbitrage potential of a Wind and Storage Power Plant (W&SPP) participating in the multistage spot energy market, and thus leading to total profit improvement [11], [12].

Many factors make decision-making for VRRGs complicated. Although sophisticated forecasting techniques exist, forecast errors are unavoidable and WFs are forced to take corrective actions, which usually result in lower profits [12], [13]. We develop a decision-making framework that considers uncertainty through the development of two decision making frameworks with stochastic optimization models.

## 2. Literature Review

Stochastic programming has been used extensively in the literature to optimize the participation of renewable, non-renewable [14], [15], and combination energy generators in electricity markets around the world. In this literature review, we will explore some scenario generation techniques as well as applications for stochastic optimization models that include schedule-related decision-making for systems that contain wind farms.

The variability for VRRG leads to uncertainty when planning energy generation. There are various methods combined with stochastic programming to deal with that uncertainty and increase the reliability of the proposed solution. The model in Ref. [16] includes chance-constrained stochastic program features in a two-stage stochastic program. In Ref. [17], robust optimization techniques are used to represent the uncertainty through confidence bounds. In Refs. [18] and [19], historical data is used to generate deterministic forecasts using Long Short-Term Memory Recurrent Neural Network (LSTM-RNN) as well as univariate and multivariate clustering-based k-means algorithms. In Ref. [20], scenario tree construction algorithms to successively reduce the number of nodes to decrease the computational burden and keep the problem tractable.

Stochastic programming models have also been used to aid decision-making from perspectives other than that of a generator. In Ref. [21], a stochastic decision-making framework is developed from the viewpoint, not of a producer, but rather a local market operator, aggregator, or prosumers. In Ref. [22], the problem is looked at from the perspective of an aggregator for electric vehicle charging stations. In Ref. [23], a two-sided two-stage optimization model simultaneously considers both the supply and the demand side of wind power to ensure stable consumption in the real-time market. Ref. [24] formulates a two-stage robust stochastic mixed-integer programming (SMIP) to obtain an optimal generation expansion plan considering generation mix and construction time to satisfy forecasted electricity demand. Ref [25] looks at the problem from the perspective of a micro-grid that receives loads from various intermittent renewable energy sources.

Due to the instability of wind power generation, many systems are a combination of wind energy with another power source. Ref. [26] develops stochastic optimal distribution scheduling models for hybrid wind-solar photovoltaic (PV) systems, wherein system uncertainties also include those that impact solar PV plants such as irradiance. In Ref. [27], instead of solar PV, the work integrates with the wind farm a concentrating solar power (CSP) plants that are equipped with thermal energy storage (TES). In Ref. [28], a coupling series of models is developing for seeking optimal hydropower and wind power strategies. Similarly, in Ref. [29], a combination of natural gas and wind power systems are considered and modeled using stochastic optimization.

One method of eliminating the negative characteristics of uncertainty for renewable energy power generation and making up for forecasting errors is by including a form of energy storage in the system [30]. Refs. [31] and [11] look at hybrid power systems consisting of WFs and batteries to co-optimize both the day-ahead offering and nominal real-time operating strategies of W&SPP. In Ref. [18], the analysis of the stochastic program includes a comparison of various sizes of ESSs such as batteries. Ref. [32] looks at an optimal bid submission in a day-ahead electricity market for the problem of joint operation of wind with photovoltaic power systems having an energy storage device. In Ref. [33], rather than a battery, hydro pumped storage units are considered in the formulation. Ref. [34] also looks at hydro pumps to manage production but does so while also including risk-hedging through Conditional Value at Risk.

It is also important to consider the participation of the generator in the reserve regulation market. In Ref. [19], an LSTM-RNN is designed to generate forecasts for regulation requirements. In Ref. [35], producers are encouraged to use the reserve market to regulate short term-output by using some of the generation mismatches as regulation reserve services instead of appearing as energy imbalance to avoid paying penalties and increase profit. Ref. [36] develops a SMIP to schedule reserves provided by demand response providers, which act as aggregators for SOs. In Refs. [37], [38] reserve dispatch is considered using a risk-averse approach by minimizing the conditional value-at-risk (CVaR), however, in Ref. [38], the model also considers the BM.

Evidently, stochastic optimization is used extensively to optimize wind farm participation in the energy market, however, there is still an opportunity for improvement. Section 3 highlights the gaps in the literature that the presented work here will address.

## 3. Contributions

The contributions of this research can be categorized into two classes: contribution towards application and contribution towards methodology development. In terms of application, many of the references listed in the literature review section of this paper optimize the participation of a generator in multiple markets [39]. For instance, the work presented in [18] optimizes the participation of a WF in the day-ahead, intra-day, and balancing market, and the work presented in [19] optimization of a WF in the day-ahead and reserve market. However, none of the references in the literature review include the day-ahead, intra-day, reserve, and balancing markets all within a single framework. As such, there is a potential for improving the participation of a renewable power generator within the power-market. Findings suggest that wind power plants that in addition to the day-ahead market are also active in the adjustment [18], real-time [40], and reserve markets [35] have the opportunity to increase profitability. Not only can the generator increase competitiveness by selling energy in the market that offers the best price, but also, the added flexibility allows the generator to participate in market-based arbitrage. This work will fill that gap by expanding on the work presented in [18] and [19], developing a model that is inclusive of all four markets.

In terms of methodology, as the time gets closer to the actual time of participation, the accuracy of the information that the decision-maker has increased, however, two-stage stochastic programming models do not typically allow for information to be updated as it is made available. Having more accurate information at the time of decision-making can greatly improve the quality of the decisions. Therefore, in this work, a second framework is developed that includes a mechanism to update scenarios and their associated probabilities through a phase-based approach. This phase-based two-stage programming approach is a novel methodological approach that allows decisions to be updated as more information becomes available to the decision-maker. Comparatively to multi-stage stochastic programming, this approach avoids the complexity of building a scenario tree in each stage while maintaining the ability to model the lookahead uncertainty and use newly available information at each decision stage.

### 4. Model Description

#### 4.1. Short-term energy market structure

In this work, we will study the operation of a wind farm (WF) in a generalized energy market. The decision-making framework follows the temporal framework imposed by the energy market. To create the most generalized framework that applies to many markets, the framework studied in this work contains all of the different forms of markets running concurrently and sequentially. We can think of the participation of W&SPP as split up between two primary categories of markets, the pool market (PM) and the reserve market, where the PM consists of the day-ahead, intra-day, and balancing market, as depicted in Fig. 1. In this model, we assume the generator assumes the role of a price-taker.

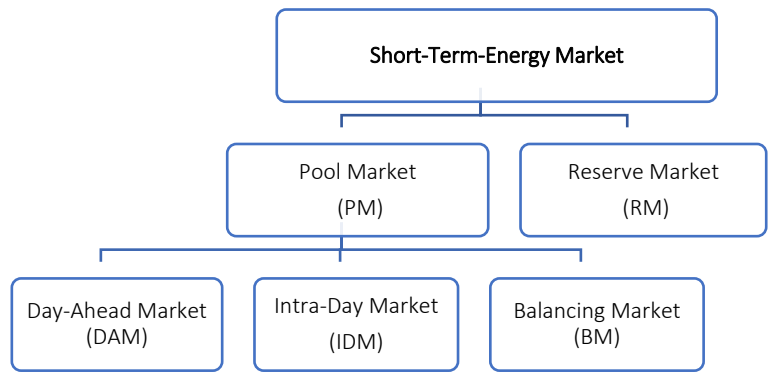

*Fig. 1.* Generalized short-term energy market structure.

The generator must commit to how much it plans to buy/sell in each market. These commitments occur sequentially leading up to the time of actual participation as shown in Fig. 2. The first commitments are the commitments made in the morning of the day before in what is called the day-ahead market (DAM). This market typically yields the best return, however, at this point, information regarding the availability of energy and market prices is inaccurate. At this point, the generator must also make commitments to the reserve market (RM) regarding how much it will allocate to fulfill spot requirements made by the SO for regulation purposes. Later, starting from the evening of the day before the intra-day market (IDM) is initiated, more accurate information regarding the market prices is available and the generator can adjust the commitments made. The generator is further able to adjust commitments throughout the next day in the balancing market (BM).

Note that some markets around the world may not contain all of these mechanisms. For instance, most European energy markets contain an adjustment market, and many American energy markets do not [8]. The work presented was made to be generalizable, so it is still applicable. In those instances, the price would be set to zero in the model.

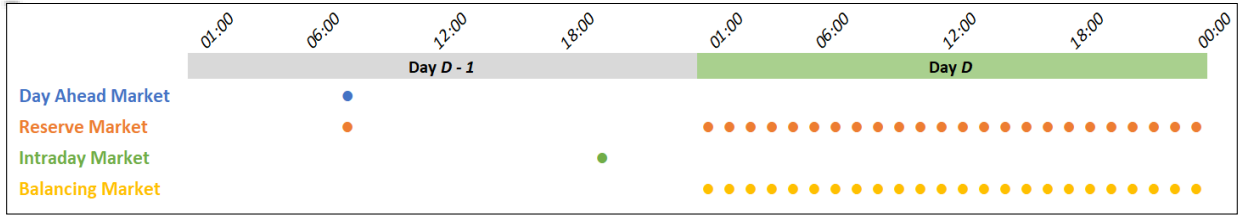

*Fig. 2. Schedule of commitments to buy/sell in various markets.*

The producer is assumed to be a price-taker, which means that there is no price interaction between the generation offer and the market-clearing price. This is because we assume that the capacity of the WF is small relative to the energy market volume and its actions do not impact the energy market outcomes [41].

In this work, several key factors stop the system from performing buying and selling actions in the same market in the same time period. For the regulation market, regulation up and down is an input. Therefore, it is the system operator that decides whether regulation up or down is needed at any point. For the balancing market, the prices are modeled to avoid arbitrage between buying and selling within the same market for the same time period. For the charging and discharging process, an efficiency parameter was included in the model. Therefore, some energy is lost during the discharging/charging, and it is never efficient to both discharge and charge in the same time period.

### 4.2. Scenario dependent programming

Stochastic programming is a popular modeling approach for problems where decisions are made in stages, and between the stages, some uncertainty in the problem parameters is unveiled. The decisions in subsequent stages depend on the outcome of the uncertain parameters. For instance, in the application of WFs, decisions about how much energy to commit must be made before knowing the amount of available wind energy (AWE), market prices, and regulation requirements will be. Then, decisions made on the actual day of participation, in real-time, will be dependent on how much was committed as well as the AWE, market prices, and regulation requirements. Thus, the problem is well suited for the stochastic programming approach. In this approach, decision variables are categorized into the first-stage and second-stage decisions as shown in Fig. 3. The first-stage decisions are those that must be made ahead of time before accurate information being available. Second-stage decisions are those that happen after the information becomes available and are dependent on that information.

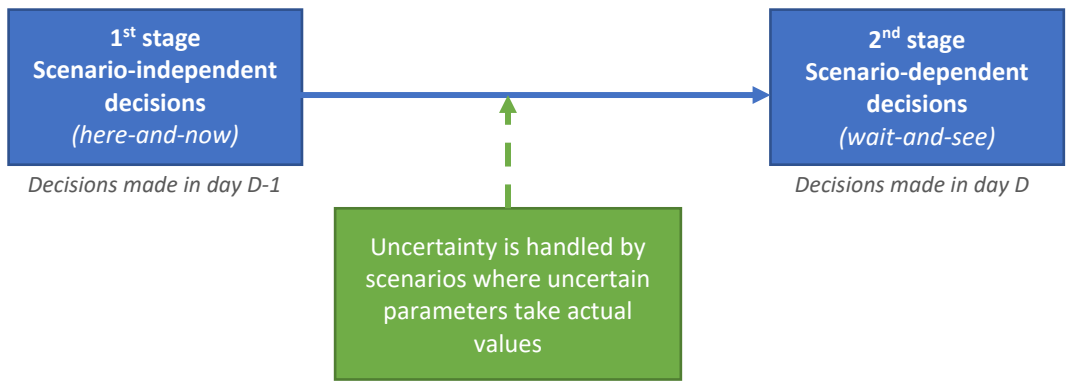

*Fig. 3. Two-stage stochastic programming overview.*

To model the uncertainty of the parameters, scenarios are developed to come up with representative predictions of what the parameter values could be, and probabilities are assigned to each of the scenarios. In this model, the uncertain parameters are the available wind energy, market prices for all markets, and the regulation requirements. It is assumed that these parameters change on an hourly basis.

### 4.3. Stochastic programming frameworks

Two frameworks will be presented in this work and the results will be compared. Framework A consists of a two-stage stochastic convex optimization approach, where both scenario-independent and scenario-dependent decisions are made in a single phase. Framework B consists of four sequential stochastic optimization models wherein the results from each model feed into each subsequent model allowing for scenarios to be updated as more information becomes available to the decision-maker. In this section, the two frameworks will be detailed.

#### 4.3.1. Framework A: Single-phase framework

The first framework is a single-phase two-stage stochastic programming model and serves as a baseline. This model is the first in literature, as far as the author's knowledge, to contain all the markets presented. All decisions made in day $D-1$, which are the decisions associated with the participation of the W&SPP in DAM, IDM, and commitments to the RM, are considered as first-stage variables as shown in Fig. 4. Participation in the BM and actual participation in RM are modeled as scenario-dependent second-stage variables.

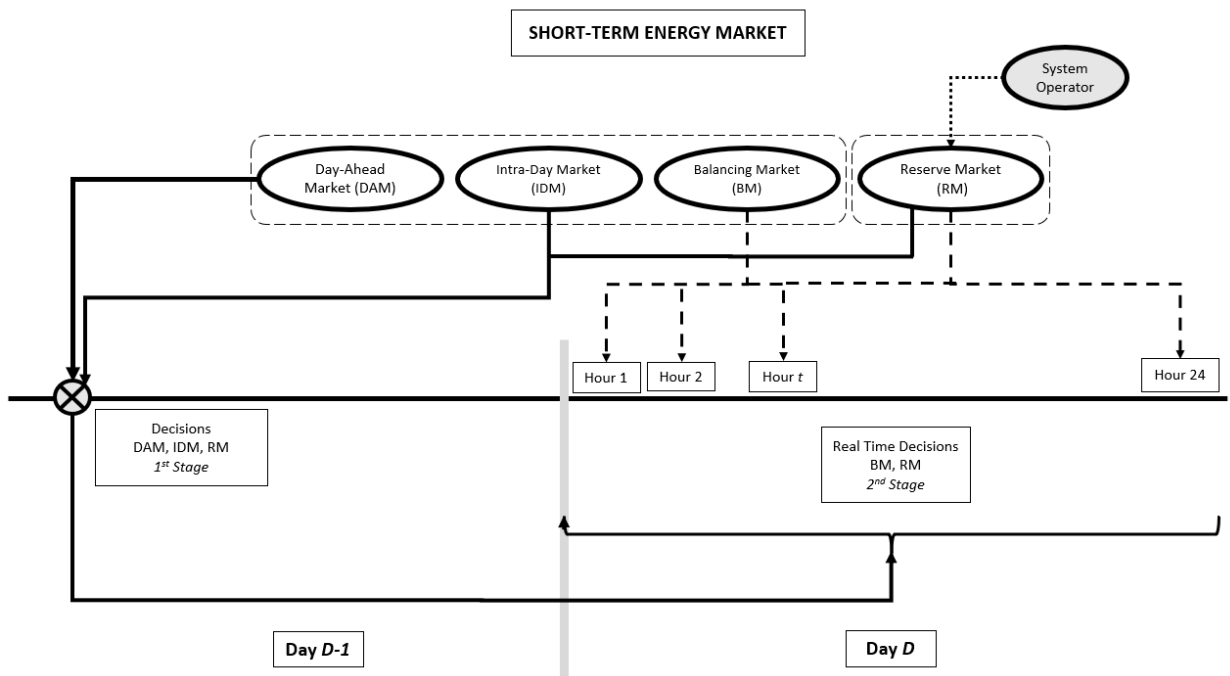

**Fig. 4.** *Baseline framework overview – Framework A.*

*4.3.2. Framework B: Multi-phase sequential framework*

One major limitation to the baseline model is that if more information becomes available, and uncertain data becomes more accurate, we are unable to use that information to update decisions. As such a multi-phase model is developed where four stochastic programming models are arranged sequentially as in Fig. 5. This approach of allowing information to be updated and decisions to be revised early in the process is novel. Outputs (decisions) from one phase are fed into the following phase, and input (uncertain) data to the problem is used as soon as it becomes available or updated.

In Fig. 5 each of the phases of Framework B is summarized, and the following sections describe each phase in detail. Note that for each phase, the objective function is to maximize profit.

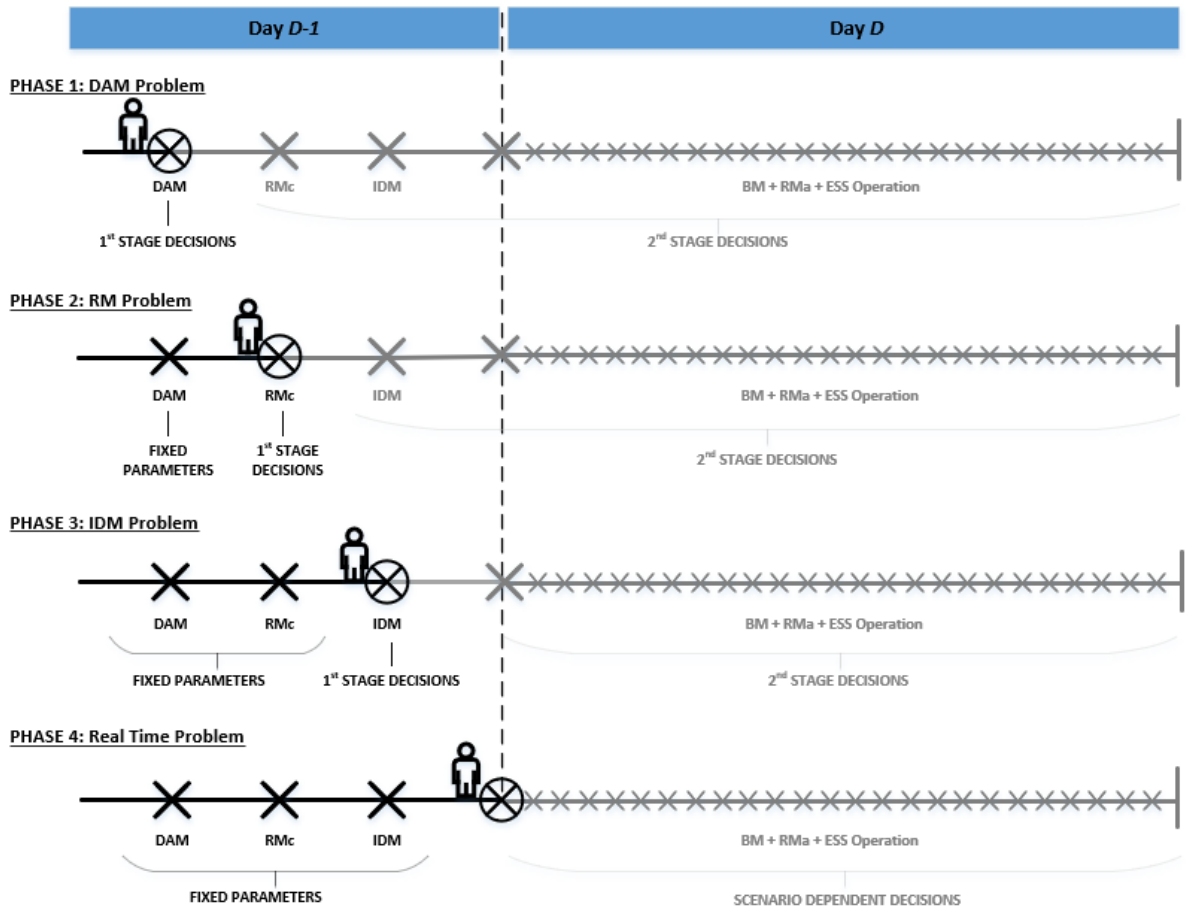

*Fig. 5. Overview of the multi-phase framework – Framework B.*

*4.3.2.1. Phase 1 – DAM Problem*

The sequence begins with the problem of deciding how much to commit in the DAM. At this stage, the available wind energy and market prices are unknown and are represented by scenarios in the model. Potential participation in other markets is represented in the model as second-stage variables.

*4.3.2.2. Phase 2 – Regulation Band Problem*

The reserve market is operated by the system operator (not the market operator as is the case for DAM and IDM). The goal of this market is the regulation of real-time power system operation. The SO imposes a regulation band, and the agent must commit to providing a certain amount of power, which may be used by the SO in real-time for regulation tasks in case it is needed. The fraction of the committed reserves that are required by the SO for regulation tasks in real-time is referred to as regulation requirements. The regulation band is a +/- amount of power that the generator may provide if required by the SO. If the SO requires regulation up, it means that the generator is supposed to increase its generation. The converse applies for regulation down. This power market also runs under a bidding mechanism and even though the system operator imposes the regulation, the system can still decide how much to participate in and whether to fulfill the required commitments or incur a penalty. In this work it is assumed that monetary penalty is incurred and is modeled by a penalty factor [42].

In this work, it is assumed that only the ESS may be used for regulation requirements. It is challenging for a wind farm to perform regulation capabilities when available wind energy is uncertain. However, the ESS allows for added operational flexibility and quick response to ancillary services such as frequency regulation [43]. Note that regulation "up" refers to selling energy to the market, and "down" refers to buying energy from the market.

*4.3.2.3. Phase 3 – IDM Problem*

The intraday market has been shown to improve wind producers' competitiveness [44] when combined with the day ahead and balancing market [45], predominantly for its role in reducing balancing needs [46].

*4.3.2.4. Phase 4 – Real-Time Problem*

The real-time problem involves the BM, as well as real-time energy offered for RM. During day $D$, for real time regulation purposes, the system operator may ask the generator to supply regulation, i.e., a percentage of the regulation band committed. This percentage is modeled by the parameter $\pi$ which is considered uncertain.

The BM prices are defined as the imbalance prices from the source. As with Framework A, we assume deviations in DAM and IDM are covered in the BM. The balancing market is a two-price settlement system where two different prices are considered for deviation up and deviation down. Similar to the regulation market, balancing market price "up" refers to selling energy, and "down" refers to buying energy.

4.4. Nomenclature

| Sets and subindices | |
|---|---|
| $S$ | The set of all scenarios under consideration |
| $T$ | The set of all periods under consideration |
| $s$ | Subindex for scenarios, $s \in S$ |
| $t$ | Subindex for time slots, $t \in T$ |

Parameters

| | |
|---|---|
| $\bar{P}^{wind}$ | Rated power of the wind farm (MW) |
| $E_0^{ess}$ | Initial energy stored in the ESS (MWh) |
| $\eta_{in}$ | Charging efficiency of the ESS |
| $\eta_{out}$ | Discharging efficiency of the ESS |
| $\bar{E}^{ess}$ | ESS maximum allowable storage capacity (MWh) |
| $\bar{P}^{ess}$ | Maximum power to/from ESS (MW) |
| $SOC^{min}$ | Minimum state of charge allowed for the ESS |
| $\rho_s$ | Probability of scenario $s$ |
| $\hat{P}_{s,t}^{wind}$ | Forecasted wind power available in every hour of scenario $s$ (MW) |
| $P_t^{prod}$ | Amount of power generated in wind farm |

Power Market Parameters

*Energy Market*

| | |
|---|---|
| $\beta_{s,t}^{dam}$ | Scenario generated energy price in the DAM in every hour of day $D$ for scenario $s$ (€/MWh) |
| $\beta_{s,t}^{idm}$ | Scenario generated energy price in the IDM in every hour of day $D$ for scenario $s$ (€/MWh) |
| $\lambda_t^{bm,up}$ | Energy price of deviation up in every hour of day $D$ (€/MWh) |
| $\lambda_t^{bm,dw}$ | Energy price of deviation down in every hour of day $D$ (€/MWh) |

*Reserve Market*

| | |
|---|---|
| $R^{rm,up}$ | Ratio between reserves up and total reserves. Constant. |
| $\kappa^{rm}$ | Correction factor for the cost of deviation in the amount offered and required for regulation up and down in RM. Constant. |
| $\gamma_{s,t}^{rm}$ | Price of power reserve in every hour of day $D$ for scenario $s$ (€/MWh) |
| $\beta_{s,t}^{rm,up}$ | Price of energy under regulation up in RM for every hour of day $D$ (€/MWh |
| $\beta_{s,t}^{rm,dw}$ | Price of energy under regulation down in RM for every hour of day $D$ (€/MWh) |
| $\lambda_{s,t}^{rm,up}$ | Cost of deviation in the amount of energy offered and required for regulation up in RM for every hour of day $D$ (€/MWh) |
| $\lambda_{s,t}^{rm,dw}$ | Cost of deviation in the amount of energy offered and required for regulation down in RM for every hour of day $D$ (€/MWh) |
| $\pi_{s,t}^{rm,up}$ | Regulation requirement up by SO in every hour of day $D$ for scenario $s$ (ratio between actual energy and reserved power) |
| $\pi_{s,t}^{rm,dw}$ | Regulation requirement down by SO in every hour of day $D$ for scenario $s$ (ratio between actual energy and reserved power) |

Decision variables

*Overall System*

| | |
|---|---|
| $P_{s,t}^{wind}$ | Amount of power produced by the wind farm (MW) |
| IDAM | Total income from participation in the day-ahead market for day D |
| IIDM | Total income from participation in the intraday market for day D |

| IBM | Total income from participation in the balancing market for day D |
| --- | --- |
| IRM | Total income from participation in the reserve market for day D |

*ESS Operation*

| $E_{s,t}^{ess}$ | Energy stored in the ESS in every hour of day *D* for scenario *s* (MWh) |
| --- | --- |
| $P_{s,t}^{ess,in}$ | Power entering to the ESS in every hour of day *D* for scenario *s* (MW) |
| $P_{s,t}^{ess,out}$ | Power delivered by the ESS in every hour of day *D* for scenario *s* (MW) |
| $SOC_{s,t}$ | State of charge of ESS in every hour of day *D* for scenario *s* |

*Energy Market*

| $\hat{P}_t^{dam}$ | Final power committed in the DAM for every hour of day *D* (MW) |
| --- | --- |
| $\hat{P}_t^{idm}$ | Final power committed in the IDM in every hour of day *D* (MW) |
| $\hat{P}_t^{pm}$ | Final power committed in the Energy Market in every hour of day *D* (MW) |
| $P_t^{pm}$ | Power actually traded in the Energy Market in every hour of day *D* (MW) |
| $\Delta_{s,t}^{bm}$ | Deviation in power committed in the Energy Market in every hour of day D for scenario s (MW) |
| $\hat{P}_{s,t}^{bm,up}$ | Final power committed for BM market up (MW) |
| $\hat{P}_{s,t}^{bm,dw}$ | Final power committed for BM market down (MW) |

*Reserve Market*

| $\hat{P}_t^{rm}$ | Final power committed for RM in every hour of day *D* (MW) |
| --- | --- |
| $\hat{P}_t^{rm,up}$ | Final power committed for regulation up in RM in every hour of day *D* (MW) |
| $\hat{P}_t^{rm,dw}$ | Final power committed for regulation down in RM in every hour of day *D* (MW) |
| $\hat{E}_{s,t}^{rm,up}$ | Energy required by SO for regulation up in RM every hour of day *D* for scenario *s* (MWh) |
| $\hat{E}_{s,t}^{rm,dw}$ | Energy required by SO for regulation down in RM every hour of day *D* for scenario *s* (MWh) |
| $E_{s,t}^{rm,up}$ | Energy offered for regulation up in RM every hour of day *D* for scenario *s* (MWh) |
| $E_{s,t}^{rm,dw}$ | Energy offered for regulation down in RM every hour of day *D* for scenario *s* (MWh) |
| $\Delta_{s,t}^{rm,up}$ | Deviation in regulation up in RM every hour of day *D* for scenario *s* (MWh) |
| $\Delta_{s,t}^{rm,dw}$ | Deviation in regulation down in RM every hour of day *D* for scenario *s* (MWh) |

4.5. Framework A

4.5.1. Objective Function of Framework A

Framework A consists of a stochastic programming model in which the objective function is to maximize the income of the operation of the system participating in all markets considered: DAM, IDM, BM, and RM.

$$maximize\ (IDAM + IIDM + IBM + IRM) \tag{1}$$

Each of the income variables in the objective function (1) is defined in equations (2)-(5). Note that the income equations are summed over the entire set of scenarios. The income for the RM (3) is a little more complex than the others in that it is a combination of the income from the amount of energy committed towards the RM, ($\sum_{t \in T} \gamma_{s,t}^{rm} \hat{P}_t^{rm}$), the income from energy actually offered for regulation in the RM, ($\sum_{t \in T} \beta_t^{rm,up} E_{s,t}^{rm,up} - \sum_t \beta_t^{rm,dw} E_{s,t}^{rm,dw}$), and the costs associated with deviating from requirements ($\sum_{t \in T} \lambda_t^{rm,up} D_{s,t}^{rm,up} - \sum_{t \in T} \lambda_t^{rm,dw} D_{s,t}^{rm,dw}$).

$$\text{IDAM} = \sum_{s \in S} \rho_s \left( \sum_{t \in T} \beta_{s,t}^{dam} \hat{P}_t^{dam} \right) \quad (2)$$

$$\text{IRM} = \sum_{s \in S} \rho_s \left( \sum_{t \in T} \gamma_{s,t}^{rm} \hat{P}_t^{rm} + \sum_{t \in T} \beta_t^{rm,up} E_{s,t}^{rm,up} - \sum_{t \in T} \beta_t^{rm,dw} E_{s,t}^{rm,dw} \right.$$
$$\left. - \sum_{t \in T} \lambda_t^{rm,up} \Delta_{s,t}^{rm,up} - \sum_{t \in T} \lambda_t^{rm,dw} \Delta_{s,t}^{rm,dw} \right) \quad (3)$$

$$\text{IIDM} = \sum_{s \in S} \rho_s \left( \sum_{t \in T} \beta_{s,t}^{idm} \hat{P}_t^{idm} \right) \quad (4)$$

$$\text{IBM} = \sum_{s \in S} \rho_s \left( \sum_{t \in T} \lambda_{s,t}^{bm,up} \Delta_{s,t}^{bm,up} - \sum_{t \in T} \lambda_{s,t}^{bm,dw} \Delta_{s,t}^{bm,dw} \right) \quad (5)$$

Substituting the income equations given in (2)-(5) into the objective function, the objective function becomes as given in equation (6).

$$maximize \sum_{s \in S} \rho_s \left( \sum_{t \in T} \beta_{s,t}^{dam} \hat{P}_t^{dam} + \sum_{t \in T} \beta_{s,t}^{idm} \hat{P}_t^{idm} + \sum_{t \in T} \gamma_{s,t}^{rm} \hat{P}_t^{rm} + \sum_{t \in T} \beta_t^{rm,up} E_{s,t}^{rm,up} \right.$$
$$- \sum_{t \in T} \beta_t^{rm,dw} E_{s,t}^{rm,dw} - \sum_{t \in T} \lambda_t^{rm,up} D_{s,t}^{rm,up} - \sum_{t \in T} \lambda_t^{rm,dw} D_{s,t}^{rm,dw}$$
$$\left. + \sum_{t \in T} \lambda_{s,t}^{bm,up} \Delta_{s,t}^{bm,up} - \sum_{t \in T} \lambda_{s,t}^{bm,dw} \Delta_{s,t}^{bm,dw} \right) \quad (6)$$

### 4.5.2. Constraints of Framework A

Operation of W&SPP:

The constraints modeling the operation of the W&SPP are the same for every phase of every framework presented in this work. Henceforth in this work, operation of W&SPP constraints will refer to equations (7)-(15).

The AWE at any given time could be more than the rated power of the farm, however, the maximum amount of power that can be produced is capped at the rated power of the farm. Therefore constraints (7) and (8) ensure that the amount of energy the WF produces, $P_{s,t}^{wind}$, is both less than total AWE, $\hat{P}_{s,t}^{wind}$, and the rated power of the farm, $\bar{P}^{wind}$.

$$P_{s,t}^{wind} \leq \hat{P}_{s,t}^{wind} \qquad \forall t \in T, \forall s \in S \qquad (7)$$
$$P_{s,t}^{wind} \leq \bar{P}^{wind} \qquad \forall t \in T, \forall s \in S \qquad (8)$$

Constraint (9) defines the amount of energy stored in ESS in every time step as a function of the initial conditions, power entering and leaving the ESS, and the efficiency of the charging and discharging processes. The efficiency of the charging, $\eta_{in}$ and discharging $\eta_{out}$ of the ESS are considered to be constants in the open interval (0,1). Note that since the charging and discharging coefficients are always considered <1 for practical purposes, then a constraint for the non-simultaneity of the charging and discharging process is not needed. This is because it is never efficient for the system to both charge and discharge the battery in the same time period.

$$E_{s,t}^{ess} = E_0^{ess} + \sum_{\tau=1}^{t} \eta_{in} P_{s,t}^{ess,in} - \sum_{\tau=1}^{t} \frac{1}{\eta_{out}} P_{s,t}^{ess,out} \qquad \forall t \in T, \forall s \in S \qquad (9)$$

Constraints (10)-(12) limit the maximum and minimum energy stored in the ESS.

$$E_t^{ess} \leq \bar{E}^{ess} \qquad \forall t \in T \qquad (10)$$
$$SOC_{s,t} = E_{s,t}^{ess} / \bar{E}^{ess} \qquad \forall t \in T, \forall s \in S \qquad (11)$$
$$SOC_{s,t} \geq SOC^{min} \qquad \forall t \in T, \forall s \in S \qquad (12)$$

Constraints (13)-(14) limit the maximum power that can be exchanged by ESS at any time.

$$P_{s,t}^{ess,out} \leq \bar{P}^{ess} \qquad \forall t \in T, \forall s \in S \qquad (13)$$
$$P_{s,t}^{ess,in} \leq \bar{P}^{ess} \qquad \forall t \in T, \forall s \in S \qquad (14)$$

Constraint (15) are non-negativity restrictions.

$$P_{s,t}^{ess,out} \,;\, P_{s,t}^{ess,in} \,;\, E_{s,t}^{ess} \geq 0 \qquad \forall t \in T, \forall s \in S \qquad (15)$$

Constraints for EM

The EM is a combination of the DAM and IDM, as in constraint (16). Constraint (17) limits the maximum power that can be bought/sold in the DAM. For the DAM, the constraints are set to respect the physical limits of the plant, meaning that the system cannot sell more than the available energy, and cannot buy more than the storage capacity. Similarly, constraint (18) limits the maximum power that can be bought/sold in the IDM. However, the IDM limits are set such that the W&SPP can completely update its position in the IDM.

$$\hat{P}_t^{pm} = \hat{P}_t^{dam} + \hat{P}_t^{idm} \qquad \forall t \in T \qquad (16)$$
$$-\bar{P}^{ess} \leq \hat{P}_t^{dam} \leq \bar{P}^{wind} + \bar{P}^{ess} \qquad \forall t \in T \qquad (17)$$
$$|\hat{P}_t^{idm}| \leq \bar{P}^{wind} + \bar{P}^{ess} \qquad \forall t \in T \qquad (18)$$

Deviations in the PM are covered by the BM as in constraint (19). The balancing market is then defined as in constraints (20) -(23).

$$\Delta_{s,t}^{bm} = P_{s,t}^{pm} - \hat{P}_{t}^{pm} \quad \forall t \in T, \forall s \in S \quad (19)$$

$$\Delta_{s,t}^{bm} = \hat{P}_{s,t}^{bm,up} - \hat{P}_{s,t}^{bm,dw} \quad \forall t \in T, \forall s \in S \quad (20)$$

$$\hat{P}_{s,t}^{bm,up} \leq \bar{P}^{wind} + \bar{P}^{ess} \quad \forall t \in T, \forall s \in S \quad (21)$$

$$\hat{P}_{s,t}^{bm,dw} \leq \bar{P}^{wind} + \bar{P}^{ess} \quad \forall t \in T, \forall s \in S \quad (22)$$

$$\hat{P}_{s,t}^{bm,up}; \hat{P}_{s,t}^{bm,dw} \geq 0 \quad \forall t \in T, \forall s \in S \quad (23)$$

Constraints for the Reserve Market:

In this work, the penalty for not fulfilling reserve market requirements is modeled by defining a parameter $\kappa^{rm} > 1$ and including the set of equations (24) and (25). In this work, it is assumed that a monetary penalty is applied in the case of non-fulfillment. In other markets, failure to supply in the reserve may result in being banned. In such cases the parameter $\kappa^{rm}$ should be assigned a sufficiently large number to ensure fulfillment of commitments.

$$\lambda_{s,t}^{rm,up} = \kappa^{rm} \beta_{s,t}^{rm,up} \quad \forall t \in T, \forall s \in S \quad (24)$$

$$\lambda_{s,t}^{rm,dw} = \kappa^{rm} \beta_{s,t}^{rm,dw} \quad \forall t \in T, \forall s \in S \quad (25)$$

Constraints (26)-(28) define the regulation band that can be offered by the W&SPP.

$$\hat{P}_{t}^{rm} = \hat{P}_{t}^{rm,up} + \hat{P}_{t}^{rm,dw} \quad \forall t \in T \quad (26)$$

$$\hat{P}_{t}^{rm,up} \leq \bar{P}^{ess} \quad \forall t \in T \quad (27)$$

$$\hat{P}_{t}^{rm,dw} \leq \bar{P}^{ess} \quad \forall t \in T \quad (28)$$

Constraint (29) defines the ratio between the regulations up and the total regulation band offered. This ratio must follow the ratio assigned for the entire system, which is a constant parameter $R$ that is given by the Spanish system operator. The rate is fixed based on the general requirement of the system. Some markets in other countries are asymmetric and do not care about the shape. In those cases, constraint (29) may be excluded.

$$\hat{P}_{t}^{rm,up} / \hat{P}_{t}^{rm} = R^{rm,up} \quad \forall t \in T \quad (29)$$

Constraints (30) and (31) set the amount of energy required by SO for regulation tasks.

$$\hat{E}_{s,t}^{rm,up} = \pi_{s,t}^{rm,up} \cdot \hat{P}_{t}^{rm,up} \quad \forall t \in T, \forall s \in S \quad (30)$$

$$\hat{E}_{s,t}^{rm,dw} = \pi_{s,t}^{rm,dw} \cdot \hat{P}_{t}^{rm,dw} \quad \forall t \in T, \forall s \in S \quad (31)$$

Constraints (32) and (33) deal with the actual energy supplied by W&SPP for regulation tasks.

$$E_{s,t}^{rm,up} \leq \hat{E}_{s,t}^{rm,up} \quad \forall t \in T, \forall s \in S \quad (32)$$

$$E_{s,t}^{rm,dw} \leq \hat{E}_{s,t}^{rm,dw} \quad \forall t \in T, \forall s \in S \quad (33)$$

Constraints (34) and (35) define deviations in the RM.

$$\Delta_{s,t}^{rm,up} = \hat{E}_{s,t}^{rm,up} - E_{s,t}^{rm,up} \quad \forall t \in T, \forall s \in S \quad (34)$$

$$\Delta_{s,t}^{rm,dw} = \hat{E}_{s,t}^{rm,dw} - E_{s,t}^{rm,dw} \qquad \forall t \in T, \forall s \in S \qquad (35)$$

Constraint (36) establishes some non-negativity requirements.

$$\hat{P}_t^{rm,up}; \hat{P}_t^{rm,dw}; E_{s,t}^{rm,up}; E_{s,t}^{rm,dw} \geq 0 \qquad \forall t \in T, \forall s \in S \qquad (36)$$

<u>Power balance for all markets</u>

The final constraint (37) ensures the power balance between the system and the power market. The left-hand-side of the equation is the total amount of energy that is being offered to the system and the right-hand-side of the equation is the total amount of energy that is available. Note that energy and power variables can be added and subtracted because an hourly time step is considered throughout the framework.

$$P_{s,t}^{pm} + E_{s,t}^{rm,up} - E_{s,t}^{rm,dw} = P_t^{wind} + P_t^{ess,out} - P_t^{ess,in} \qquad \forall t \in T, \forall s \in S \qquad (37)$$

This concludes Framework A. In the following section, Framework B is given.

### 4.6. Framework B: Four sequential phases

Framework B consists of four sequential stochastic programming models that feed into one another. Phase one is the DAM problem, phase two is the RM problem, phase three is the IDM problem, and phase four is the real-time problem. One way of summarizing Framework B is by looking at which variables are scenario dependent, scenario independent, and which variables become fixed parameters as outlined in Table 1 and visualized in Fig. 5.

After each phase is solved, the resulting first-stage decisions are fed into the subsequent problem in the form of fixed parameters. Then scenarios are updated using more up-to-date forecasts. One of the second stage decision variables is made into a 1st stage decision variable and the model is run again.

| | Decision Variables | | |
|---|---|---|---|
| | **Fixed from the previous phase** | **Scenario Independent** 1st Stage | **Scenario Dependent** 2nd Stage |
| **Phase 1:** (Unknown market prices, Unknown AWE) | - | $\hat{P}_t^{dam}$ | $\hat{P}_{s,t}^{idm}, \Delta_{s,t}^{bm},$ $\hat{P}_{s,t}^{rm}, \hat{E}_{s,t}^{rm}$ |
| **Phase 2:** (DAM prices known, other market prices unknown, Updated AWE) | $\hat{P}_t^{dam}$ | $\hat{P}_t^{rm}$ | $\hat{P}_{s,t}^{idm}, \Delta_{s,t}^{bm},$ $\hat{E}_{s,t}^{rm}$ |
| **Phase 3:** (DAM prices known, RM price is known, other market prices unknown, Updated AWE) | $\hat{P}_t^{dam},$ $\hat{P}_t^{rm}$ | $\hat{P}_t^{idm}$ | $\Delta_{s,t}^{bm},$ $\hat{E}_{s,t}^{rm}$ |
| **Phase 4:** (All market prices known, Updated AWE forecasts) | $\hat{P}_t^{dam}, \hat{P}_t^{idm}$ $\hat{P}_t^{rm}$ | - | $\Delta_{s,t}^{bm},$ $\hat{E}_{s,t}^{rm}$ |

*Table 1. Scenario-dependent and scenario-independent variables for the multi-phase model.*

### 4.6.1. Phase 1 of Framework B

Phase 1 is the DAM problem.

1st stage variables: commitments made in DAM in day $D-1$

2nd stage variables: participation in RM, IDM, BM, ESS operation.

Note, at this stage all prices are unknown, and forecasts are inaccurate.

$$maximize\ (\text{IDAM} + \text{IIDM} + \text{IBM} + \text{IRM}) \qquad (38)$$

$$\text{IDAM} = \sum_{s \in S} \rho_s \left( \sum_{t \in T} \beta_{s,t}^{dam} \hat{P}_t^{dam} \right) \qquad (39)$$

$$\text{IRM} = \sum_{s \in S} \rho_s \left( \sum_{t \in T} \gamma_{s,t}^{rm} \hat{P}_{s,t}^{rm} + \sum_{t \in T} \beta_t^{rm,up} E_{s,t}^{rm,up} - \sum_{t \in T} \beta_t^{rm,dw} E_{s,t}^{rm,dw} \right. \qquad (40)$$
$$\left. - \sum_{t \in T} \lambda_t^{rm,up} \Delta_{s,t}^{rm,up} - \sum_{t \in T} \lambda_t^{rm,dw} \Delta_{s,t}^{rm,dw} \right)$$

$$\text{IIDM} = \sum_{s \in S} \rho_s \left( \sum_{t \in T} \beta_{s,t}^{idm} \hat{P}_{s,t}^{idm} \right) \qquad (41)$$

$$\text{IBM} = \sum_{s \in S} \rho_s \left( \sum_{t \in T} \lambda_{s,t}^{bm,up} \Delta_{s,t}^{bm,up} - \sum_{t \in T} \lambda_{s,t}^{bm,dw} \Delta_{s,t}^{bm,dw} \right) \qquad (42)$$

Substituting the Phase 1 income equations (39)-(42) into the objective function (38), the objective function becomes as given in (43).

$$maximize \sum_{s \in S} \rho_s \left( \sum_{t \in T} \beta_{s,t}^{dam} \hat{P}_t^{dam} + \sum_{t \in T} \gamma_{s,t}^{rm} \hat{P}_{s,t}^{rm} + \sum_{t \in T} \beta_t^{rm,up} E_{s,t}^{rm,up} \right. \qquad (43)$$
$$- \sum_{t \in T} \beta_t^{rm,dw} E_{s,t}^{rm,dw} - \sum_{t \in T} \lambda_t^{rm,up} \Delta_{s,t}^{rm,up} - \sum_{t \in T} \lambda_t^{rm,dw} \Delta_{s,t}^{rm,dw}$$
$$\left. + \sum_{t \in T} \beta_{s,t}^{idm} \hat{P}_{s,t}^{idm} + \sum_{t \in T} \lambda_{s,t}^{bm,up} \Delta_{s,t}^{bm,up} - \sum_{t \in T} \lambda_{s,t}^{bm,dw} \Delta_{s,t}^{bm,dw} \right)$$

W&SPP operational constraints are as described in Framework A, given by equations (7)-(15).

### 4.6.2. Phase 2 of Framework B

RM Problem. At this point, decisions have already been made as to how much the system will commit to the DAM, and there is new information regarding the requirements made by the system operator. AWE forecast is updated. DAM prices are known.

Again, the objective is to maximize total income.

Decisions regarding DAM have already been made and will not be modified during this phase.

1st stage variables: commitments to RM made in day $D-1$

2nd stage variables: participation in IDM, BM, ESS operation, and deviations in participation in RM.

$$maximize\ (\text{IIDM} + \text{IBM} + \text{IRM}) \tag{44}$$

$$\text{IRM} = \sum_{s \in S} \rho_s \left( \sum_{t \in T} \gamma_{s,t}^{rm} \hat{P}_t^{rm} + \sum_{t \in T} \beta_t^{rm,up} E_{s,t}^{rm,up} - \sum_{t \in T} \beta_t^{rm,dw} E_{s,t}^{rm,dw} \right.$$
$$\left. - \sum_{t \in T} \lambda_t^{rm,up} D_{s,t}^{rm,up} - \sum_{t \in T} \lambda_t^{rm,dw} D_{s,t}^{rm,dw} \right) \tag{45}$$

$$\text{IIDM} = \sum_{s \in S} \rho_s \left( \sum_{t \in T} \beta_{s,t}^{idm} \hat{P}_{s,t}^{idm} \right) \tag{46}$$

$$\text{IBM} = \sum_{s \in S} \rho_s \left( \sum_{t \in T} \lambda_{s,t}^{bm,up} \Delta_{s,t}^{bm,up} - \sum_{t \in T} \lambda_{s,t}^{bm,dw} \Delta_{s,t}^{bm,dw} \right) \tag{47}$$

Substituting the Phase 2 income equations given in (45)-(47) into the objective function (44), the objective function becomes as given in (48).

$$maximize \sum_{s \in S} \rho_s \left( \sum_{t \in T} \gamma_{s,t}^{rm} \hat{P}_t^{rm} + \sum_{t \in T} \beta_t^{rm,up} E_{s,t}^{rm,up} - \sum_{t \in T} \beta_t^{rm,dw} E_{s,t}^{rm,dw} \right.$$
$$- \sum_{t \in T} \lambda_t^{rm,up} \Delta_{s,t}^{rm,up} - \sum_{t \in T} \lambda_t^{rm,dw} \Delta_{s,t}^{rm,dw} + \sum_{t \in T} \beta_{s,t}^{idm} \hat{P}_{s,t}^{idm}$$
$$\left. + \sum_{t \in T} \lambda_{s,t}^{bm,up} \Delta_{s,t}^{bm,up} - \sum_{t \in T} \lambda_{s,t}^{bm,dw} \Delta_{s,t}^{bm,dw} \right) \tag{48}$$

W&SPP operational constraints are as described in Framework A, given by equations (7)-(15).

### 4.6.3. Phase 3 of Framework B

IDM Problem

Decisions regarding the DAM and RM have already been made and will not be modified here.

1st stage variables: commitments to IDM made in day $D-1$

2nd stage variables: participation in BM, ESS operation, and deviations in participation in IDM.

$$maximize\ (\text{IIDM} + \text{IBM} + \text{IRM}) \tag{49}$$

$$\text{IRM} = \sum_{s \in S} \rho_s \left( \sum_{t \in T} \beta_t^{rm,up} E_{s,t}^{rm,up} - \sum_{t \in T} \beta_t^{rm,dw} E_{s,t}^{rm,dw} - \sum_{t \in T} \lambda_t^{rm,up} D_{s,t}^{rm,up} \right.$$
$$\left. - \sum_{t \in T} \lambda_t^{rm,dw} D_{s,t}^{rm,dw} \right) \tag{50}$$

$$\text{IIDM} = \sum_{s \in S} \rho_s \left( \sum_{t \in T} \beta_t^{idm} \hat{P}_t^{idm} \right) \tag{51}$$

$$\text{IBM} = \sum_{s \in S} \rho_s \left( \sum_{t \in T} \lambda_{s,t}^{bm,up} \Delta_{s,t}^{bm,up} - \sum_{t \in T} \lambda_{s,t}^{bm,dw} \Delta_{s,t}^{bm,dw} \right) \tag{52}$$

Substituting the Phase 3 income equations (50)-(52) into the objective function (49), the objective function becomes equation (53).

$$maximize \sum_{s \in S} \rho_s \left( \sum_{t \in T} \beta_t^{rm,up} E_{s,t}^{rm,up} - \sum_{t \in T} \beta_t^{rm,dw} E_{s,t}^{rm,dw} - \sum_{t \in T} \lambda_t^{rm,up} \Delta_{s,t}^{rm,up} \right.$$
$$- \sum_{t \in T} \lambda_t^{rm,dw} \Delta_{s,t}^{rm,dw} + \sum_{t \in T} \beta_t^{idm} \hat{P}_{s,t}^{idm} + \sum_{t \in T} \lambda_{s,t}^{bm,up} \Delta_{s,t}^{bm,up}$$
$$\left. - \sum_{t \in T} \lambda_{s,t}^{bm,dw} \Delta_{s,t}^{bm,dw} \right) \tag{53}$$

W&SPP operational constraints are as detailed in Framework A, given by equations (7)-(15).

### 4.6.4. Phase 4 of Framework B

Real-time problem – note, there are no scenario-independent variables.

$$maximize \text{ (IBM + IRM)} \tag{54}$$

$$\text{IRM} = \sum_{s \in S} \rho_s \left( \sum_{t \in T} \beta_t^{rm,up} E_{s,t}^{rm,up} - \sum_{t \in T} \beta_t^{rm,dw} E_{s,t}^{rm,dw} - \sum_{t \in T} \lambda_t^{rm,up} D_{s,t}^{rm,up} \right.$$
$$\left. - \sum_{t \in T} \lambda_t^{rm,dw} D_{s,t}^{rm,dw} \right) \tag{55}$$

$$\text{IBM} = \sum_{s \in S} \rho_s \left( \sum_{t \in T} \lambda_{s,t}^{bm,up} \Delta_{s,t}^{bm,up} - \sum_{t \in T} \lambda_{s,t}^{bm,dw} \Delta_{s,t}^{bm,dw} \right) \tag{56}$$

Substituting the Phase 4 income equations (55)-(56) into the objective function (54), the objective function becomes equation (57).

$$maximize \sum_{s \in S} \rho_s \left( \sum_{t \in T} \beta_t^{rm,up} E_{s,t}^{rm,up} - \sum_{t \in T} \beta_t^{rm,dw} E_{s,t}^{rm,dw} - \sum_{t \in T} \lambda_t^{rm,up} \Delta_{s,t}^{rm,up} \right.$$
$$\left. - \sum_{t \in T} \lambda_t^{rm,dw} \Delta_{s,t}^{rm,dw} + \sum_{t \in T} \lambda_{s,t}^{bm,up} \Delta_{s,t}^{bm,up} - \sum_{t \in T} \lambda_{s,t}^{bm,dw} \Delta_{s,t}^{bm,dw} \right) \tag{57}$$

W&SPP operational constraints are as detailed in Framework A, given by equations (7)-(15).

This concludes Framework B. In the following section, the simulation procedure is described, and some results are shared along with some insights from the results.

## 5. Results

The goal of this section is to describe and evaluate the two proposed frameworks. To do this, a real-world example is used. A 60-day evaluation is conducted for a Spanish WF starting on a random day within a 1-year dataset of prices, wind energy, and regulation requirements. Market prices and SO regulation requirements are obtained from a Spanish System Operator Information System [47] that is a part of the Iberian Electricity Market. Wind energy forecast was provided by the WF operator of Sotavento in Northwestern Spain [48].

Scenarios are generated for each of the three unknown parameters and are fed into the stochastic model. The stochastic model produces a theoretical expected net income given the optimal decisions generated. To evaluate the quality of the decisions generated from the stochastic model, the decisions are fed into a deterministic model using the actual values from historical data of the parameters for each day that was run.

In the following sections, the scenarios and their generation methods are described, followed by some results and key insights.

### 5.1. Scenario generation to handle uncertainty

Stochastic programming was selected as the modeling approach for this work due to the uncertainty associated with the values of certain parameters at the time of the decision-making, i.e., a day ahead of the actual participation of a wind farm in the energy market. Stochastic programming allows for the inclusion of a set of scenarios that aim to represent the possible values of the parameters with a certain probability.

There are a myriad of scenario generation approaches that are used in combination with stochastic programming. In references [18] and [19], various scenario generation approaches are implemented and compared. Comparing scenario generation approaches is not the objective of this work, hence one approach is used for generating AWE scenarios, and another one for generating scenarios with market prices and regulation requirements. The specific approaches are discussed in the following sections. It is to be noted, however, the flexibility of the framework in that it is capable of handling scenarios generated from different methodologies and of different sizes.

For simplicity, it is assumed that the unknown parameters are independent of each other. Thus, the total number of scenarios is a product of the number of scenarios for each respective parameter. In this work, the scenario tree contains 3x 10 x 3 = 90 scenarios.

#### 5.1.1. Wind Energy Scenarios

One year of historic probabilistic forecasts for available wind energy was provided by the WF operator as a set of time series on a percentile basis. Several forecasts for each day are made throughout the day prior. For framework A, the forecast that is produced closest to the time of participation is used. For Framework B, we select four forecasts, one for each phase, with equal intervals between each update. Each forecast consists of a time series corresponding to a set of hourly values defining an upper bound on the actual available wind energy with a given probability as illustrated in Fig. 6. Three scenarios are generated from curves p75, p50 and, p25 for each selected forecast for each day. Thus, the hourly AWE scenarios for each day that is run in the model are specific to that day.

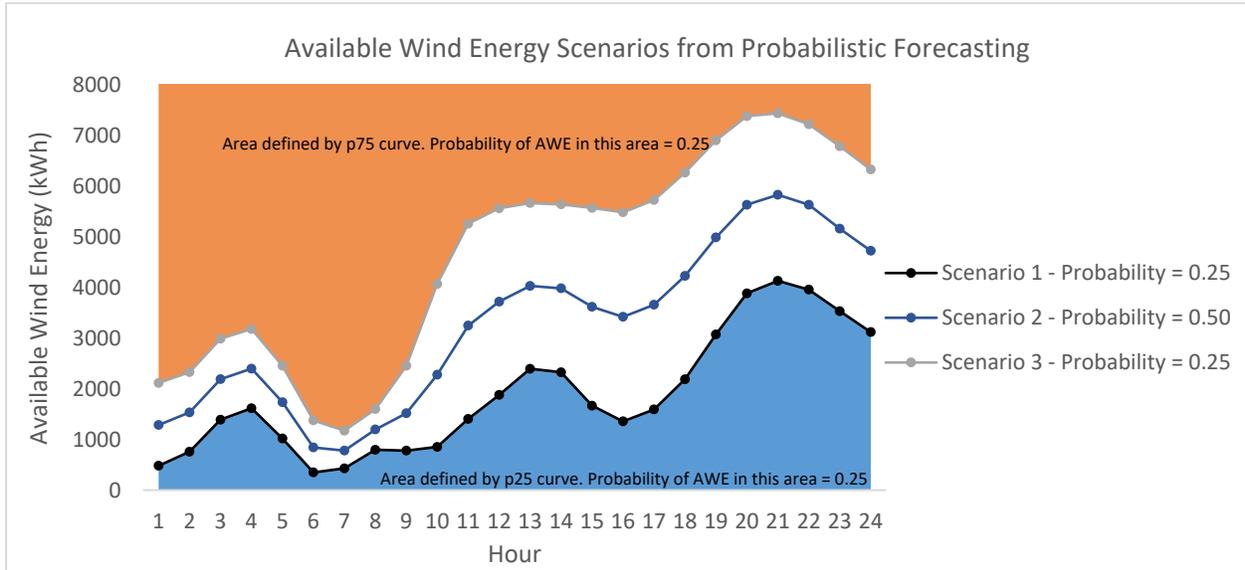

*Fig. 6. Scenarios of available wind energy from probabilistic forecasts.*

### 5.1.2. Market Price Scenarios

For market prices, we combine in one scenario a full day of data for each market, therefore one scenario comprises 24 hours x 7 market prices = 168 attributes. We do this by taking one year of historical data and arranging them into 365 vectors, where each vector is 1 day with 168 attributes. A k-means algorithm separates the vectors into 10 clusters. The centroid of each cluster is used as scenarios, and probabilities are assigned to the scenarios using frequentist reasoning, i.e., with the probability of each scenario set to equate the number of vectors assigned to each cluster as a ratio to the total number of vectors. Detailed results are available in Ref. [18].

### 5.1.3. Regulation Requirement Scenarios

The k-means algorithm is also used for the regulation requirements, but we take a simpler approach where each data point contains two attributes, regulation up and regulation down, and is then separated into 3 clusters. The scenarios are set to be the centroids of each cluster, and the probability of the occurrence of each scenario is based on the number of observations that fall in each cluster. Detailed results are available in Ref. [19].

### 5.2. Decision analysis framework comparison

A quick way to evaluate the overall performance of the two frameworks is to compare the resulting net income from the decisions produced. Fig. 7 compares the net income from a random series of 60 days for a single WF. We see that when running the same data through the two frameworks, Framework B performs better than Framework A in 87% of the runs. The mean net income increase is 635.77 Euro per day, which is equivalent to a 7% improvement.

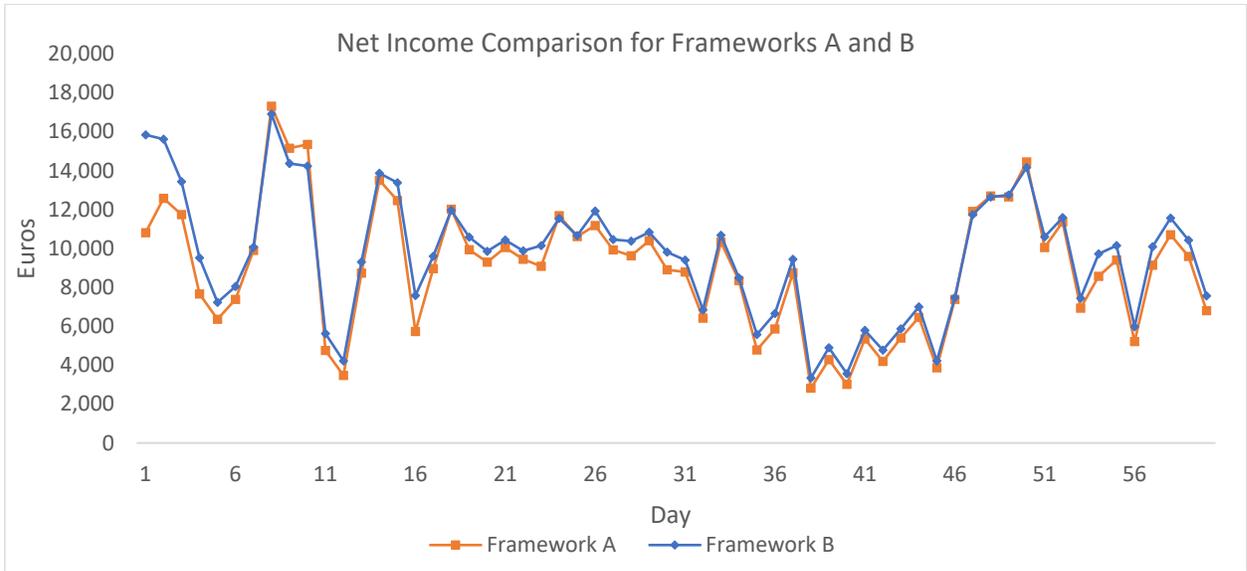

*Fig. 7.* Result comparison of Framework A and Framework B for 60 simulated days for one WF.

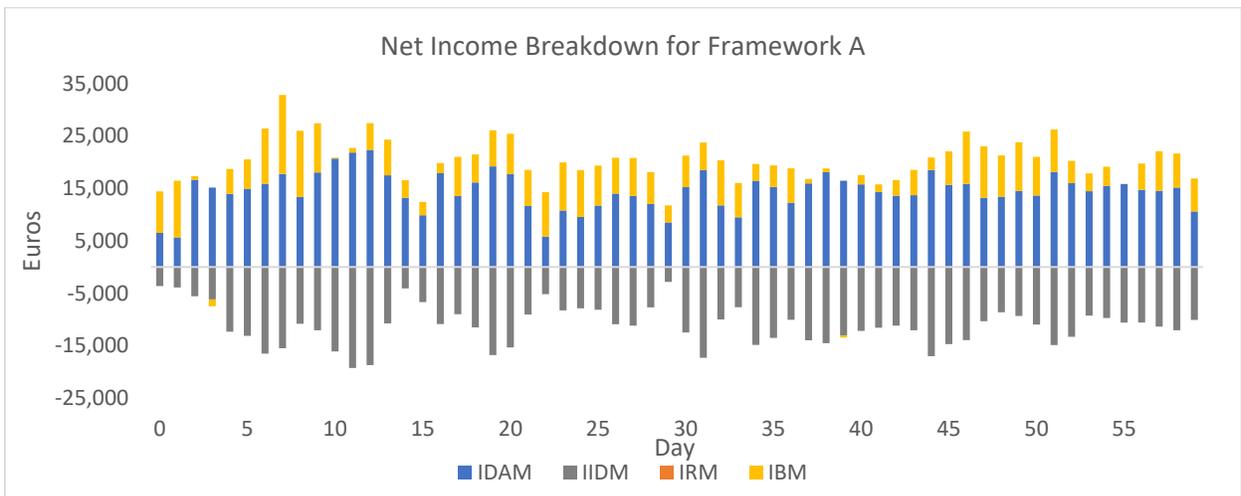

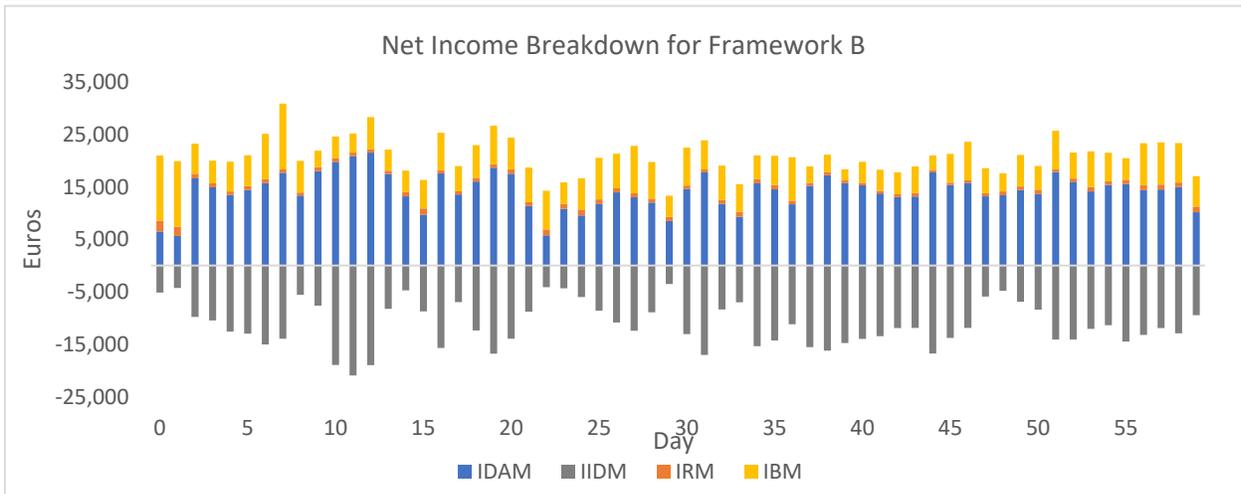

*Fig. 8.* Proportion of overall participation of the wind farm in each market for Frameworks A and B.

To understand the difference in income, Fig. 8 presents the contribution of the participation in each market to the total income for each framework, and Fig. 9 compares the income from each market for each framework side-by-side. One point that becomes immediately apparent is the increase in participation in RM in Framework B, whereas in Framework A participation in the RM is almost zero. Another intriguing point is that the IDM contributes negatively to the total income in every run for both frameworks. To explain this phenomenon, Fig. 10 presents the operational decisions for one random day *D* within the 60-day evaluation. The actual available wind energy and market prices for the given day are shown in the two upper-subfigures. The third and fourth subfigures represent the participation of the wind farm in the various markets given the decisions that are produced using frameworks A and B, respectively.

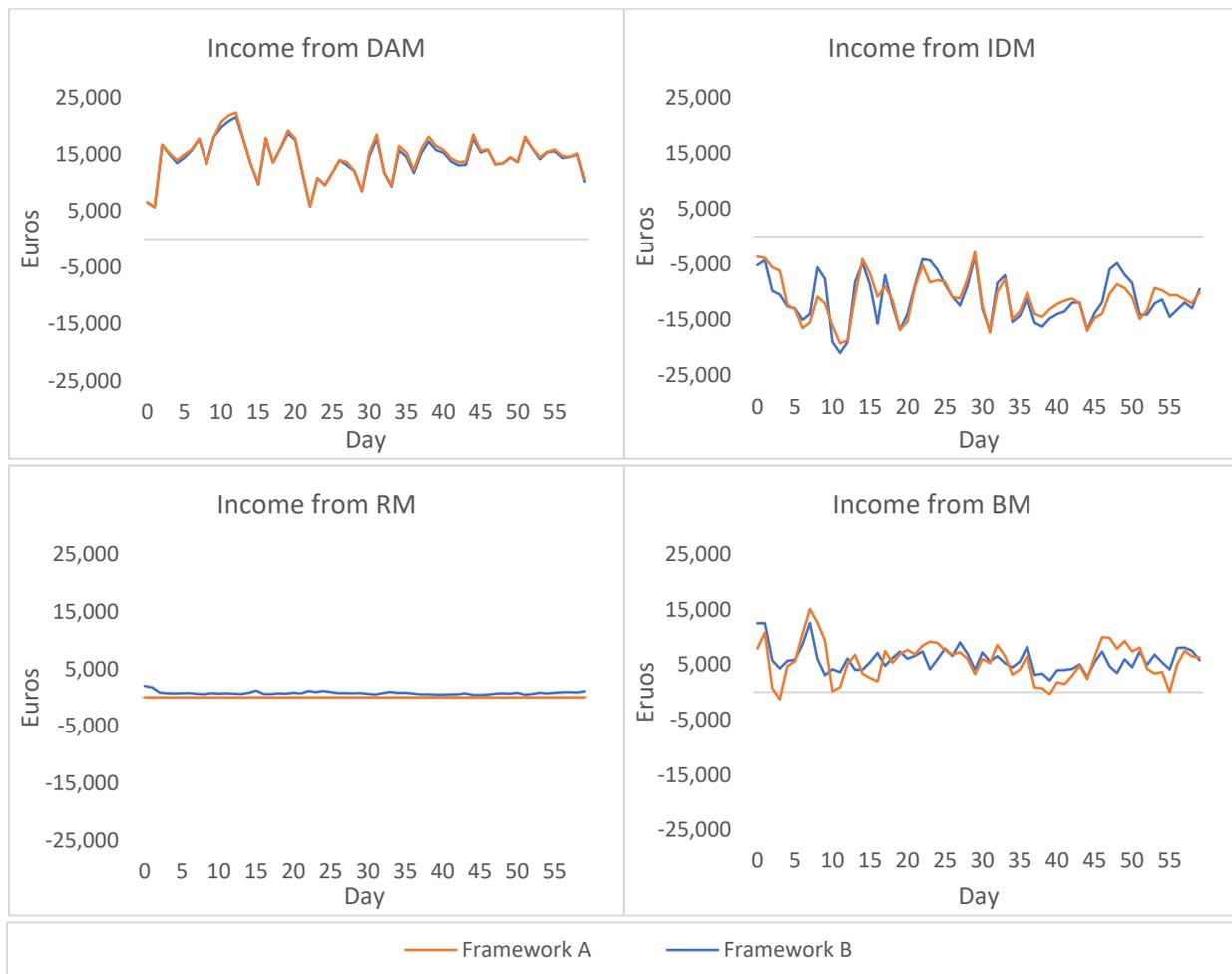

*Fig. 9. Overall participation of the wind farm in each market for Frameworks A and B.*

Although the graphs presented show the results for one day, some interesting dynamics can be seen. Further analysis shows that these dynamics are consistent in the other simulated days. For instance, in both frameworks the wind farm prioritizes selling in the DAM even at a cost of purchasing energy from the IDM to make up for the production deficiency due to unavailable wind energy. It does so to the extent that the farm sells the maximum allowable energy it can, which is capped at the rated power of the wind farm. Although the participation patterns look similar for the two frameworks, one difference is in the lower participation in the BM in Framework B. This is because at the time decision making occurs for the IDM, more information is available. The increased accuracy results in smaller deviations and less costs incurred in the BM.

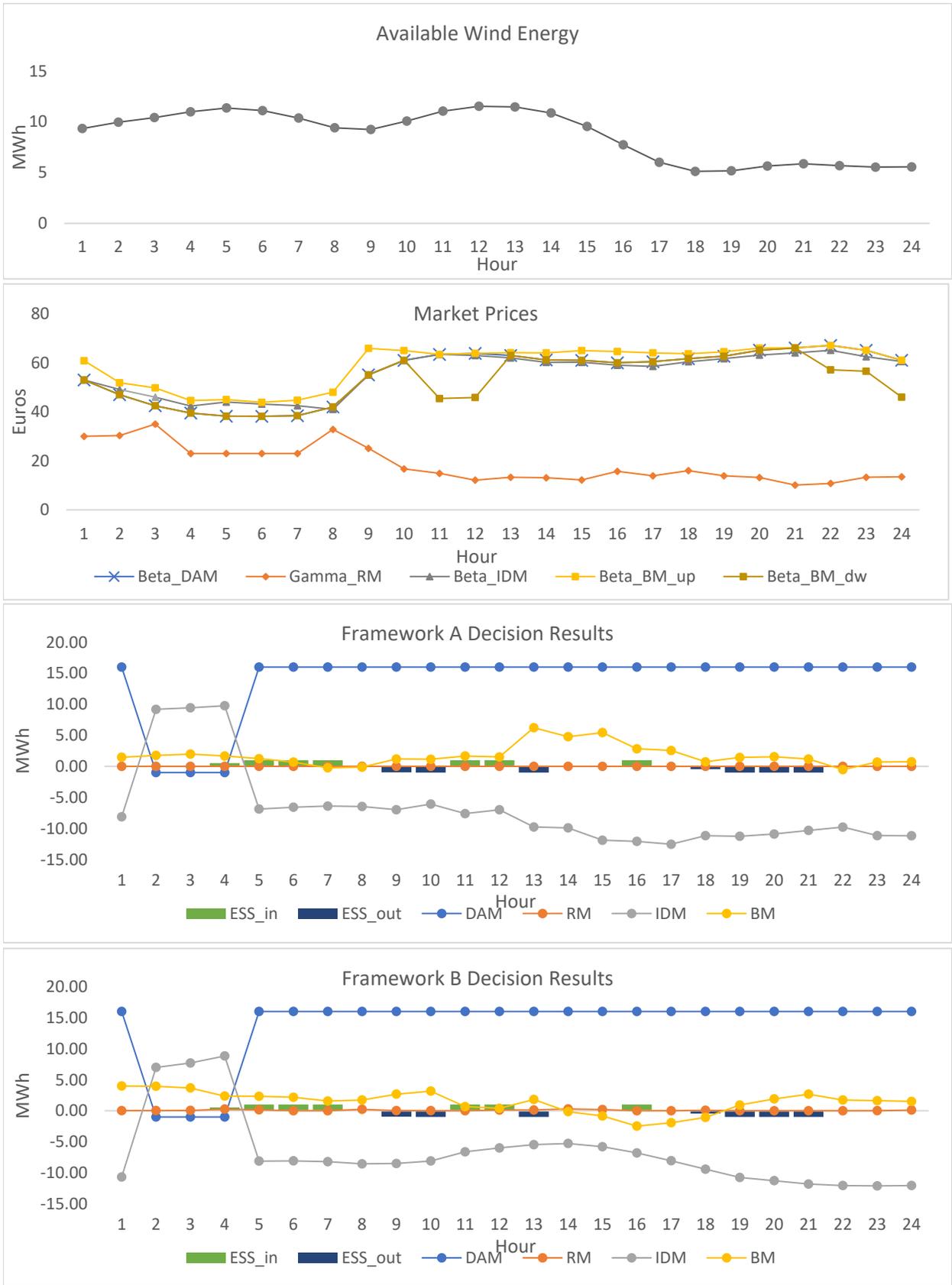

Fig. 10. Results of one day's operation of a single wind farm.

## 5.3. Computational Effort

The improvement in results from Framework A compared to Framework B comes at a cost of increased computational effort. In Framework B, multiple models must be run, one for each phase, resulting in average total CPU time approximately four times the CPU time required for Framework A. The average CPU time to run one day for Framework A and Framework B respectively are 106 seconds and 427 seconds. However, in a practical setting, the four phases of Framework B would not be run at the same time. Each sequential phase would be run as more information is made available to update decisions that have not had to be confirmed yet. The average CPU time for each phase of Framework B is summarized in Table 2. The time decreases for each sequential phase because fewer second-stage decision variables remain from one phase to the next.

|         | Average CPU time for one day (s) |
|---------|----------------------------------|
| Phase 1 | 125                              |
| Phase 2 | 110                              |
| Phase 3 | 98                               |
| Phase 4 | 93                               |

*Table 2. CPU time for each phase of Framework B.*

## 6. Conclusion

The work presented is aimed at expanding upon frameworks created to optimize the participation of wind energy producers in multiple energy markets. It did so by developing two decision-making frameworks using stochastic optimization models. The first consisted of a single model in which decisions regarding four different markets were made concurrently while accounting for unknown parameters through probability-assigned scenarios. The second framework consisted of four stochastic optimization models, in which scenario-independent decisions made in one model fed into successive models, and scenarios were updated as more information became available to the decision-maker.

Using the phase-based methodology presented in this paper, decisions can be updated as more information becomes available. The results presented in this paper demonstrate the higher quality of decisions provided when using the phase-based methodology. In the experiment that we ran, we saw a mean improvement of 7% in the profitability of the farm.

The comparison between the two frameworks demonstrated the importance of updated information as forecasts become more accurate when the time comes closer to actual participation. The main observation is that since the penalty for overgeneration is less than the spot market price, the producer commits to the maximum amount that it can and makes up for any shortfalls by buying in the other markets regardless of the wind forecast.

The frameworks also demonstrated the increased economic benefit of allowing the energy producer to participate in all the markets. In each simulation run, the participation of the wind energy producer in each one of the markets was significant to the strategy, decisions made, and overall profit.

The frameworks were created in a generalizable way so they can be applied to most electricity markets around the world. They are also applicable to energy sources other than wind farms, such as solar PV, since they follow the same market dynamics.

There are several limitations to these frameworks that could translate to opportunities for future work. Firstly, the multi-phase framework required more computational effort when compared to the single-phase framework. Further improvement to the model could be made by attempting to reduce the computational load. This could be achieved by reducing the number of scenarios. One way of doing this without omitting important information is to question the assumption of independent scenarios. In reality, the unknown parameters could be dependent, such as market prices being driven by the availability of energy. Thus, scenarios could be created that include information regarding all unknown parameters. The size of each scenario will increase, but the number of scenarios would decrease, leading to a smaller number of variables and constraints.

Secondly, in this work, the producer is assumed to be a price-taker, which is appropriate for modeling a single or a small number of WF. However, if the producer is large enough, bids could influence the market. Thus, an extension to the current model would be to relax the price-taker constraint. Another opportunity for improvement would be to perform a case study on forecasting methods of wind power, regulation reserve, and prices data, and apply the different methods to the frameworks presented. Lastly, the proposed models exclude the amount of risk associated with the proposed solution. Future work could quantify and minimize the risk of the solution offered.

**Acknowledgments**

The authors would like to thank Victoria Okoria Tonbara, an industrial engineering graduate student at the department of Mechanical and Industrial Engineering, Northeastern University for her help in proofreading the manuscript. The authors would also like to thank the editor and the anonymous reviewers for their constructive comments to improve the overall quality of this paper.